\documentclass[12pt]{amsart}

\usepackage{amssymb}

\theoremstyle{plain}

\newtheorem{thm}{Theorem}[section]

\newtheorem{lmm}[thm]{Lemma}

\theoremstyle{definition}

\newtheorem{rmk}[thm]{Remark}

\def\pmc#1{\setbox0=\hbox{#1}
    \kern-.1em\copy0\kern-\wd0
    \kern.1em\copy0\kern-\wd0}

\def\de{\delta}

\def\si{\sigma}

\def\vp{\varphi}

\def\om{\omega}

\def\op{\operatorname}

\def\ov{\overline}

\def\sm{\setminus}

\author{Saharon Shelah}

\address{Department of Mathematics, 
Hebrew University, Givat Ram, 91904 Jerusalem, Israel}

\author{Katsuya Eda}
\address{School of Science and Engineering, 
Waseda University, Tokyo 169-0072, Japan}
\email{eda@logic.info.waseda.ac.jp}

\begin{document}

\title{The non-commutative Specker phenomenon in the uncountable case}

\maketitle

\section{}
An infinitary version of the notion of free products has been 
introduced and investigated by G. Higman \cite{Higman:unrestrict}. Let
$G_{i} (i\in I)$ be groups and $\ast_{i\in X} G_{i}$ the free product of
$G_{i} (i\in X)$ for $X \Subset I$ and $ p _{XY}: \ast_{i\in Y} G_{i}
\rightarrow \ast_{i\in X} G_{i}$ the canonical homomorphism for
$X\subseteq Y \Subset I.$ ( $X\Subset I$ denotes that $X$ is a finite
subset of $I$.)  

Then, the unrestricted free product is the inverse limit $\varprojlim
(\ast_{i\in X} G_{i},  p_{XY}: X \subseteq Y \Subset I).$  

We remark $\ast_{i\in \emptyset} G_{i} = \{ e\}$.

For the simplicity, we abreviate $\varprojlim(\ast_{i\in X} G_{i},
p_{XY}: X\subseteq Y \Subset I)$ and $\varprojlim(\ast_{n\in \om} 
{\mathbb Z}_n,  p_{mn}: m\le n<\om)$ by $\varprojlim \ast G_i$ and
$\varprojlim \ast {\mathbb Z}_n$ in the sequel. 

For $S\subseteq I$, $p_S : \varprojlim \ast G_i \to \varprojlim \ast
G_i$ be the canonical projection defined by:  $p_S(x)(X) = x(X\cap S)$
for $X\Subset I$.

\begin{thm}\label{thm:main}
Let $F$ be a free group. Then, for each homomorphism $h: 
\varprojlim \ast G_i \to F$ there exist countably complete ultrafilters
 $u_0,\cdots ,u_m$ on $I$ such that $h = h\cdot p_{U_0\cup \cdots 
\cup U_m}$ for every $U_0\in u_0,\cdots ,U_m\in u_m$. 

If the cardinality of the index set $I$ is less than the least
 measurable cardinal, then there exists a finite subset $X_0$ of $I$ and a
 homomorphism $\ov{h}:\ast _{i\in X_0}G_i\to F$ such that 
$h = \ov{h}\cdot p_{X_0}$, where $p_{X_0}:\varprojlim \ast G_{i}
\to \ast _{i\in X_0}G_i$ is the canonical projection. 	
\end{thm}

Previously, the second author showed the failure of the Specker
phenomenon in the uncountable case in a different situation
\cite{ShelahStruengmann:specker}. (See also \cite{ConnerEda:free}.)
We explain the difference between this result and Theorem~\ref{thm:main} 
of the present paper. There is a canonical subgroup of the unrestricted
free product, which is called the free complete product and denoted by
$\pmc{$\times$}\, \, \; _{i\in I}G_i$. When an index set $I$ is
countable, according to the Higman theorem (Lemma~\ref{lmm:higman} and a 
variant for $\pmc{$\times$}\, \, \; _{n<\om}{\mathbb Z}_n$
\cite[p.80]{Higman:unrestrict}), a homomorphism from $\varprojlim \ast
G_i $ or $\pmc{$\times$}\, \, \; _{i\in I}G_i$ to a free group factors
through a finite free product $\ast _{i\in F}G_i$. 
On the other hand, when the index set $I$ is uncountable and each $G_i$
is non-trivial, there exists a free retract of $\pmc{$\times$}\, \, \;
_{i\in I}G_i$ of large cardinality and there are homomorphisms not
factoring through any finite free product $\ast _{i\in F}G_i$, which
contrasts with the case when $I$ is countable. This also contrasts with
an abelian case, which is known as the \L o\' s theorem
\cite{EklofMekler:almostfree}. 
Theorem~\ref{thm:main} says that differing from the case of the free
complete product the non-commutative Specker phenomenon holds for the
unrestricted free products similarly as in the abelian case.  

Since the following lemma holds for the free $\si$-product 
$\pmc{$\times$}\, \, \; ^{\si}_{i\in I}{\mathbb Z}_i$ instead of a free
group $F$ \cite{E:n-specker}, Theorem~\ref{thm:main} also holds for it. 
(We remark $\pmc{$\times$}\, \, \; ^{\si}_{i\in I}{\mathbb Z}_i = 
\pmc{$\times$}\, \, \; _{i\in I}{\mathbb Z}_i$, when $I$ is countable.) 
\begin{lmm}\label{lmm:higman}$($G. Higman \cite{Higman:unrestrict}$)$
For each homomorphism $h: \varprojlim \ast {\mathbb Z}_n \to F$ there
 exists $m<\om$ and a homomorphism $\ov{h}:\ast _{i<m}{\mathbb Z}_i\to F$ 
such that $h = \ov{h}\cdot p_m$, where
 $p_m:\varprojlim \ast {\mathbb Z}_n\to \ast _{n<m}{\mathbb Z}_n$ is the
 canonical projection. 
\end{lmm}

\begin{lmm}\label{lmm:joint}
Let $I = \bigcup \{ I_n: n<\om \}$ with $I_n \subseteq I_{n+1}$ and 
$x_n \in G$ be such that $p_{I_n}(x_n) = e$. 
Then, there exists a homomorphism 
$\vp :\varprojlim \ast {\mathbb Z}_n\to \varprojlim \ast G_i$ 
such that $\vp (\de _n) = x_n$ for each $n<\om$. 
\end{lmm}

For a homomorphism $h: \varprojlim \ast G_i \to F$, let 
$\op{supp}(h)= \{ X\subseteq I: p_X(g) = e\mbox{ implies }h(g)= 0 \mbox{ for
each }g\}$. In the sequel we assume that $h$ is non-trivial. We remark
the following facts:
\begin{enumerate}
\item $p_X\cdot p_Y = p_{X\cap Y}$ for $X,Y \subseteq I$;  
\item $\op{supp}(h) = \{ X\subseteq I: h(g) = h(p_X(g))
      \mbox{ for each }g\}$;  
\item $\op{supp}$ is a filter on $I$.
\end{enumerate}

\begin{lmm}\label{lmm:complete}
Let $A_n\subseteq A_{n+1}\subseteq I$ and $A= \bigcup \{ A_n:n<\om\}$
 and $B_{n+1}\subseteq B_n\subseteq I$ and $B= \bigcap \{ B_n:n<\om\}$. 
If $A_n\notin \op{supp}(h)$ for each $n$, then $A\notin \op{supp}(h)$
 and if $B_n\in \op{supp}(h)$ for each $n$, then $B\in \op{supp}(h)$.
\end{lmm}

\begin{proof}
Suppose that $A\in \op{supp}(h)$. Take $g_n$ so that $h(g_n)\neq 0$ and
 $p_{A_n}(g_n) = e$ for each $n$ and let $u_n = p_A(g_n)$. Since 
$I = \bigcup\{ A_n\cup (I \sm A): n<\om\}$ and 
$p_{A_n\cup (I \sm A)}(u_n) = p_{A_n}(g_n) = e$, by
 Lemma~\ref{lmm:joint} we have a homomorphism  
$\vp :\varprojlim \ast {\mathbb Z}_n\to \varprojlim \ast G_i$ such that
 $\vp (\de _n) = u_n$ for each $n<\om$. Then, $h\cdot \vp (\de _n)\neq 0$ 
for each $n$, which contradicts Lemma~\ref{lmm:higman}. 

To show the second proposition by contradiction, suppose that 
$B\notin \op{supp}(h)$. Then, we have $g\in G$ such that $p_B(g) = e$
but $h(g)\neq 0$. Let $v_n = p_{B_n}(g)$. Since 
$I = \bigcup\{ B\cup (I \sm B_n): n<\om\}$ and 
$p_{B\cup (I \sm B_n)}(v_n) = p_B\cdot p_{B_n}(g_n) = e$, we apply
Lemma~\ref{lmm:joint} and have a homomorphism 
$\vp :\varprojlim \ast {\mathbb Z}_n\to \varprojlim \ast G_i$ such that
 $\vp (\de _n) = v_n$ for each $n<\om$. Then, we have a contradiction
 similarly as the above. 
\end{proof}

\begin{lmm}\label{lmm:critical}
Let $A_0\notin \op{supp}(h)$. Then, there exist $A$ satisfying the following: 
\begin{enumerate}
\item $A_0\subseteq A\notin \op{supp}(h)$; 
\item for $X\subseteq I$, $A\cup X\notin \op{supp}(h)$
      imply $(I\sm X)\cup A \in \op{supp}(h)$. 
\end{enumerate}
\end{lmm}

\begin{proof}
We construct $A_n\notin \op{supp}(h)$ by induction as follows. 
Suppose that we have constructed $A_n\notin \op{supp}(h)$. 
If $A_n$ satisfies the required properties of $A$, we
 have finished the proof. 
Otherwise, there exist $A_n\subseteq A_{n+1}\subseteq I$ 
such that $A_{n+1}\notin \op{supp}(h)$ and $(I\sm A_{n+1})\cup A_n 
\notin \op{supp}(h)$.  
We claim that this process finishes in a finite step. 
Suppose that the process does not stop in a finite step. Then, we have 
$A_n$'s and so let $A= \bigcup \{ A_n: n<\om \}$. Then, 
$A\notin \op{supp}(h)$ by Lemma~\ref{lmm:complete}.  
Since $I\sm A\subseteq I\sm A_{n+1}$, $(I\sm A)\cup A_n\notin \op{supp}(h)$ 
for each $n<\om$. Now, $I = \bigcup \{ (I\sm A)\cup A_n: n<\om \}$ and 
and hence $I\notin \op{supp}(h)$ by Lemma~\ref{lmm:complete}, which is a
 contradiction. 
\end{proof}

{\it Proof of\/} Theorem~\ref{thm:main}. 

Let $h: \varprojlim \ast G_i \to F$ be a non-trivial homomorphism. Apply
Lemma~\ref{lmm:critical} for $A_0=\emptyset$ and we have $A$. 
We define $u_0$ as follows:

$X\in u_0$ if and only if $A\cup X\in \op{supp}(h)$ for $X\subseteq
I$. Then, $u_0$ is a countably complete ultrafilter
on $I$ by Lemma~\ref{lmm:complete}. We let $I_0 = I\sm A$, then 
obviously $I\sm I_0\notin \op{supp}(h)$. 

When $I_0\in \op{supp}(h)$, then $h = h\cdot p_{U_0}$ for any $U_0\in
u_0$ and we have finished the proof. Otherewise, we construct $I_n\notin
\op{supp}(h)$ and countably complete ultrafilters 
$u_n$ on $I$ with $I_n\in u_n$ by induction as follows. 
Suppose that $\bigcup _{i=0}^nI_i\notin \op{supp}(h)$, we aplly 
Lemma~\ref{lmm:critical} for $A_0 = \bigcup _{i=0}^nI_i \notin
\op{supp}(h)$ and get a countably complete ultrafilter $u_{n+1}$ on $I$
wit $I_{n+1}\in u_{n+1}$ so that $I\sm I_{n+1}\notin \op{supp}(h)$. 

To show that this procedure stops in a finite step, suppose the
negation. Since $(I\sm \bigcup _{k=0}^\infty I_k)\cup  \bigcup _{k=0}^n
I_k$ is disjoint from $I_{n+1}$, $(I\sm \bigcup _{k=0}^\infty I_k)\cup 
\bigcup _{k=0}^n I_k\notin \op{supp}(h)$ for each $n$. Then, we have 
$I\notin \op{supp}(h)$ by Lemma~\ref{lmm:complete}, which is a
contradiction.

Now, we have a finite partition $I_0,\cdots ,I_n$ of $I$. 
By the construction, $X\in u_k$ if and only if $\bigcup _{i\neq
k}I_i\cup X\in \op{supp}(h)$ and hence for $U_k\in u_k$ ($0\le k\le n$),
then $\bigcup _{i\neq k}I_i\cup U_k\in \op{supp}(h)$. Since
$\op{supp}(h)$ is a filter, $\bigcup _{k =0}^n U_k\in \op{supp}(h)$ and
we have the first proposition. 

If each $u_j$ contains a singleton $\{ i_j\}$, we have $\{ i_0, \cdots
,i_k\}\in  
\op{supp}(h)$ and consequently the second one.
\qed

\begin{rmk}
(1) As we explained it, in the uncountable case the unrestricted free
 product and the free complete product behave differently concerning the 
 Specker phenomenon. Therefore, there is a part in the above proof which 
 cannot be converted to the case of the free complete product. It is an
 application of Lemma~\ref{lmm:joint}. The corresponding lemma to
 Lemma~\ref{lmm:joint} is \cite[Proposition 1.9]{E:free}, which holds
 under a more strict condition than Lemma~\ref{lmm:joint}, and we cannot 
 apply it. 

\noindent
(2) In \cite[theorem 1.2]{E:n-specker}, we treated with general inverse
 limits for a countable index set. For general limits we cannot
 generalize to the uncountable case \cite[Remark 1.2]{E:n-specker}. 
\end{rmk}

\providecommand{\bysame}{\leavevmode\hbox to3em{\hrulefill}\thinspace}

\end{document}